\documentclass[12pt]{article}
\usepackage{setspace}

\usepackage{amsmath}
\usepackage{amssymb}
\usepackage{multirow}
\usepackage{color}
\renewcommand{\epsilon}{\varepsilon}

\usepackage{comment}

\newcommand{\X}{\textbf{X}}

\usepackage{amsmath}
\usepackage{amssymb}
\usepackage[ansinew]{inputenc}
\usepackage{graphicx}
\usepackage{epsfig}

\newcommand{\XX}{\mathcal{X}}

\newcommand{\be}{\begin{eqnarray}}
\newcommand{\ee}{\end{eqnarray}}
\newcommand{\bea}{\begin{eqnarray*}}
\newcommand{\eea}{\end{eqnarray*}}


\title{Kernels and RKHS. Classification}

\date{}
\begin{document}

\title{Covering of high-dimensional sets }

\author{Anatoly Zhigljavsky \and Jack Noonan }

\maketitle


 \\ \

\subsection*{Introduction}
Let $(\XX,\rho)$ be a metric space and $\lambda$ be a Borel measure on this space defined on the $\sigma$-algebra generated by open subsets of $\XX$; this measure $\lambda$ defines volumes of Borel subsets of $\XX$.  The principal  case is where $\XX = \mathbb{R}^d$, $\rho $ is the Euclidean metric, and $\lambda$ is the Lebesgue measure. In this article, we are not going to pay much attention to the case of small dimensions $d$ as the problem of construction of good covering schemes for small $d$ can be attacked by the brute-force optimization algorithms. On the contrary, for medium or large dimensions (say, $d\geq 10$), there is little chance of getting  anything sensible  without understanding the main issues related to construction of efficient covering designs.

\subsection*{Optimal covering }

Let  $\X $ be a compact  subset of $\XX$  with   $0< {\rm vol}(\X)<\infty$; in order to avoid unnecessary technical difficulties, we assume that $\X$ is convex.
Consider    $X_n=\{x_1, \ldots, x_n\} $,  a set of $n$ points in $\XX$; we will
call $X_n$ an $n$-point design. The number of points $n$ can either be fixed or determined in the course of computations. In the latter case, the designs $X_n$ are incremental (nested).
The covering radius of $\X$ for the design $X_n$ is
\be \label{eq:CR}
{\rm CR} (X_n) := \max_{x\in\X}   \rho (x,X_n) \, ,
\ee
where
\be \label{eq:CR5}
\rho (x,X_n)= \min_{x_j\in X_n} \rho (x,x_j)\,
\ee
is the Hausdorff distance between  point $x \in \XX$ and  the design $X_n$.
Covering radius is also
the smallest  $r \geq 0$ such that the union of the balls with centers at $x_j \in X_n$ and radius~$r$ fully covers $\X$:
\be \label{eq:CR1}
{\rm CR} (X_n)= \mbox{$\min_{ {r>0} }$} \; \mbox{ such that }  \X \subseteq {\cal B} (X_n,r)\, ,
\ee
where ${\cal B} (X_n,r)= \bigcup_{j=1}^n {\cal B} (x_j,r)$ and
\be \label{eq:CR2} {\cal B} (x,{ r })= \{ z \in \XX :\; \rho (x,z) \leq { r } \}\ee
is the ball of radius $r$ and centre $x\in \XX$.

{\rm CR}-optimal design $X_n^{({\rm CR})}$
 is the $n$-point design such that  \bea {\rm CR}(X_n^{({\rm CR})})=\min_{X_n}  {\rm CR}(X_n). \eea

Other common names for the covering radius are:  fill distance (in approximation theory; see \cite{schaback2006kernel,wendland2004scattered}),  dispersion (in Quasi Monte Carlo; see \cite[Ch. 6]{Niederreiter}),
  minimax-distance criterion (in computer experiments;    see \cite{pronzato2012design,santner2003design}) and
 coverage threshold (in probability theory; see \cite{penrose2021random}).

Point sets with small covering radius are very desirable in theory and practice  of global optimization and many branches of numerical mathematics.
In particular, the celebrated results of A.G.Sukharev  imply that any {\rm CR}-optimal design $X_n^{({\rm CR})}$  provides the following:
(a)
 min-max $n$-point global optimization method in the set of all adaptive $n$-point optimization strategies, see \cite{Sukh1} and \cite[Ch.4,Th.2.1]{sukharev2012minimax},
 (b) worst-case  $n$-point multi-objective global optimization method in the set of all adaptive $n$-point algorithms, see
 \cite{vzilinskas2013worst}, and
  (c) the $n$-point min-max optimal quadrature, see \cite[Ch.3,Th.1.1]{sukharev2012minimax}.
  In all three cases, the class of (objective) functions is the class of Liptshitz functions, and
    the optimality of the design is independent of the value of the Liptshitz constant. Sukharev's results on $n$-point min-max optimal quadrature formulas  have been generalized in \cite{pages1998space}
    for functional classes different from the class of Liptshitz functions; see also  formula (2.3) in  \cite{du1999centroidal}.

If $\X$ is compact then choosing points outside $\X$ cannot improve best covering designs, see Proposition 3.2.3 in \cite{borodachov2019discrete} for the case $\XX=\mathbb{R}^d$, and therefore without loss of generality we can assume $x_j \in \X$ for all $j$.

In  case $\X=[0,1]$ and Euclidean metric, the $n$-point design $X_n^{({\rm CR})}= \{ x_1^*, \ldots, x_n^*\}$ minimizing ${\rm CR} (\X_n)$ consists of points $x_j^*=(2j-1)/(2n-1)$; $j=1, \ldots, n$.
The asymptotically {\rm CR}-optimal sequence of nested designs $\X_n$ for $\X=[0,1]$ can be constructed with the so-called Ruzsa points (see \cite[p.~154]{Niederreiter}).
In the case when $\X$ is a $d$-dimensional sphere,   {\rm CR}-optimal $n$-point designs (for specific values of $n$) and algorithms of  numerical construction of $n$-point designs minimizing ${\rm CR} (\X_n)$ (for the Euclidean metric) are provided  in \cite{borodachov2019discrete}.

For $d=2$, $\X=[0,1]^2$ and certain small values of $n$, either {\rm CR}-optimal or close to {\rm CR}-optimal designs are given in \cite{brass2005research,boroczky2004finite}. For $d\geq4$ and $n=4$, optimal constructions can be found in \cite{kuperbergball}.

Numerical construction of $n$-point {\rm CR}-optimal designs $X_n^{({\rm CR})}$ is notoriously difficult, especially when $d$ is not too small and $\X $ has  boundary (for example, $\X$ is neither a sphere nor a torus). This is related to the complexity of the problem of minimizing \eqref{eq:CR} with respect to $X_n$. The optimization problem is a minmax problem in a very high-dimensional  space $\XX^n$.

 Rather than constructing {\rm CR}-optimal designs, it is practically easier to regularize the optimization problem by replacing the criterion \eqref{eq:CR} with an easier one
 in such a way that the solution of a regularized problem stays close to the solution of the original problem. Both $\min$ in \eqref{eq:CR} and $\max$ in \eqref{eq:CR5} can be regularized; for example, by means of approximating the $L_\infty$-norm by a suitable $L_p$-norm.
  For examples of applications of this approach, see  \cite{borodachov2019discrete} and \cite{pronzato2019measures}.

 The problems of quantization and weak covering considered below are the two problems with two very  natural relaxations of the {\rm CR} criterion.

\subsection*{Quantization and approximate covering }

The problem of covering is a particular instance of the problem of space-filling and other space-filling criteria such as various discrepancies, separation (or packing) radius and spacing radius can be considered. However, most of these criteria (unless it is a specially constructed discrepancy) cannot be considered as regularized covering radius and therefore the designs, optimal to such criteria, can be very poor with respect to the {\rm CR}-criterion; see e.g. Section 1.3 in \cite{zhigljavsky2021bayesian} and \cite{pronzato2020bayesian}.

In what follows, it is convenient to use the intersection of the ball  \eqref{eq:CR2} and $\X$:
\be \label{eq:CR3} {B} (x,{ r })={\cal B} (x,{ r })\cap \X= \{ z \in \X :\; \rho (x,z) \leq { r } \}\, .\ee
Intersection of $n$ balls ${\cal B} (x_j,r)$ ($x_j \in X_n$) with $\X$ is therefore
\bea
{B} (X_n,r)=
{\cal B} (X_n,r) \cap \X=  \bigcup_{x_j \in X_n} { B} (x_j,r) = \{ x \in \X:\; \rho (x, X_n) \leq r\}\, ,
\eea
 where  $\rho (x,X_n) $ is  defined in \eqref{eq:CR5}.

The following cdf (cumulative distribution function) is of prime importance in understanding covering properties of the design $X_n$:
\be \label{eq:CR7}
F (r,X_n)=\mbox{vol}( {B} (X_n,r))/ \mbox{vol}( \X)\, \;\; (0 \leq r \leq {\rm CR}(X_n))\, .
\ee
For given $r\geq 0$, $F (r,X_n)$ is
the proportion of $\X$ covered by the balls ${\cal B} (X_n,r)$. For any compact $\X$ and any point set $X_n$, the distribution with cdf $F (\cdot,X_n)$ is absolutely continuous on    0 $[0,{\rm CR}(X_n)]$ so that
$F (0,X_n)=0$, $F ({\rm CR}(X_n),X_n)=1$ and the cdf $F (r,X_n)$ itself is a strictly increasing continuous function on $[0,{\rm CR}(X_n)]$.

Let $\xi=\xi (X_n)$ be the r.v. (random variable) with cdf  $F (r,X_n)$ defined by \eqref{eq:CR7}. The essential supremum ${\rm ess\, sup} \xi$ of $\xi$ is the covering radius ${\rm CR}(X_n)$. For any real $p>0$, the $p$-th moment  of $\xi$ is the so-called {\it quantization error} of order $p$:
\be
\label{eq:errorQ}
\theta_p(X_n)=\int_0^{{\rm CR}(X_n)} r^p d F (r,X_n) = \mathbb{E}_U \rho ^p(U,X_n)\, ,
\ee
where $U$ is the r.v. with uniform distribution $P_{un}$ on $\X$; that is, for any Borel subset $A$ of $\X$,
\bea
{\rm Prob}\{U \in A\}={\rm vol}(A)/{\rm vol}(\X)=P_{un}(A)\,.
\eea

As the cdf $F (r,X_n)$ itself is  strictly increasing  on $[0,{\rm CR}(X_n)]$ and ${\rm ess\, sup} \xi= {\rm CR}(X_n)$, $\theta_p(X_n)^{1/p} \to   {\rm CR}(X_n)$ for any point set $X_n$; for details, see Section 10 in \cite{graf2007foundations}. Therefore,
$p$-th order quantization error with large $p$ can be considered as a regularized covering radius.

In general, quantization does not have to be associated with uniform distribution $P_{un}$; any other probability distribution can be used in its place.
Quantization error is a very important concept with long a rich history and very important practical implications.

Quantization error is easier to numerically optimize than  covering radius. There is, for example, the celebrated `Lloyd's algorithm' which is one of the main tools in computational data science and mostly used for clustering.
For details on the theory of quantization and construction of efficient quantizers, we refer to  the excellent book \cite{graf2007foundations}.

Another important concept related to covering is the concept of {\it weak (or approximate) covering} introduced in \cite{us,second_paper} and defined through quantiles of the cdf  \eqref{eq:CR7} as follows. For  any $\gamma \in (0,1)$ and a radius $r=r_{1-\gamma}>0$, we say that  the union of  $n$ balls ${ B}(X_n,r)$ makes a
 $(1-\gamma)$-covering of $\X$ if
$
F (r,X_n)=1-\gamma \,.
$
Complete (full) covering corresponds to $\gamma=0$. As the cdf $F (\cdot,X_n)$ is continuous, $r_{1-\gamma} \to {\rm CR}(X_n)$ as $\gamma\to 0$ and therefore the problem of $(1-\gamma)$-covering with small $\gamma$ can also be considered as a regularized version of the problem of full covering. Furthermore, numerical checking of weak covering (with an approximate value of $\gamma$) is straightforward while numerical checking of the full covering is practically impossible, if $d$ is large enough. 

\subsection*{Covering of high dimensional sets }

We will now demonstrate three interesting phenomena of covering high dimensional sets; these properties are consequences of the papers \cite{us,second_paper,noonan2022efficient}.  The first phenomena is that asymptotic properties (as $n\rightarrow \infty$) are extremely far from being reached in a (realistic) finite $n$ regime. Consequently, the asymptotic results produce poor approximations in high dimensions even for large $n$ like $n=100,000$. Let us demonstrate this now. For $X_n$ a collection of $n$ i.i.d. uniform points in $\X$, the asymptotic behaviour of $F (\cdot,X_n)$ as $n\rightarrow \infty$  was the focus of study in \cite{us}.
In particular, we have the following asymptotic result:
\be
\label{eq:Zador1}
 F_{n,d}(t) \rightarrow F_d(t):= 1-\exp(-t^{d})    \;\; \mbox{ as $n\rightarrow \infty$}\,,
\ee
where $F_{n,d}(t):={\rm Pr}( n^{1/d} V_d^{1/d} \rho (U,X_n) \leq t ) $ and $V_d= {\pi}^{d/2} / \left[\Gamma (d/2\!+\!1)\right]\,$
 is the volume of the unit Euclidean ball ${\cal B}(0,1)$.

\begin{figure}[!h]
\centering
\begin{minipage}{.5\textwidth}
  \includegraphics[width=1\linewidth]{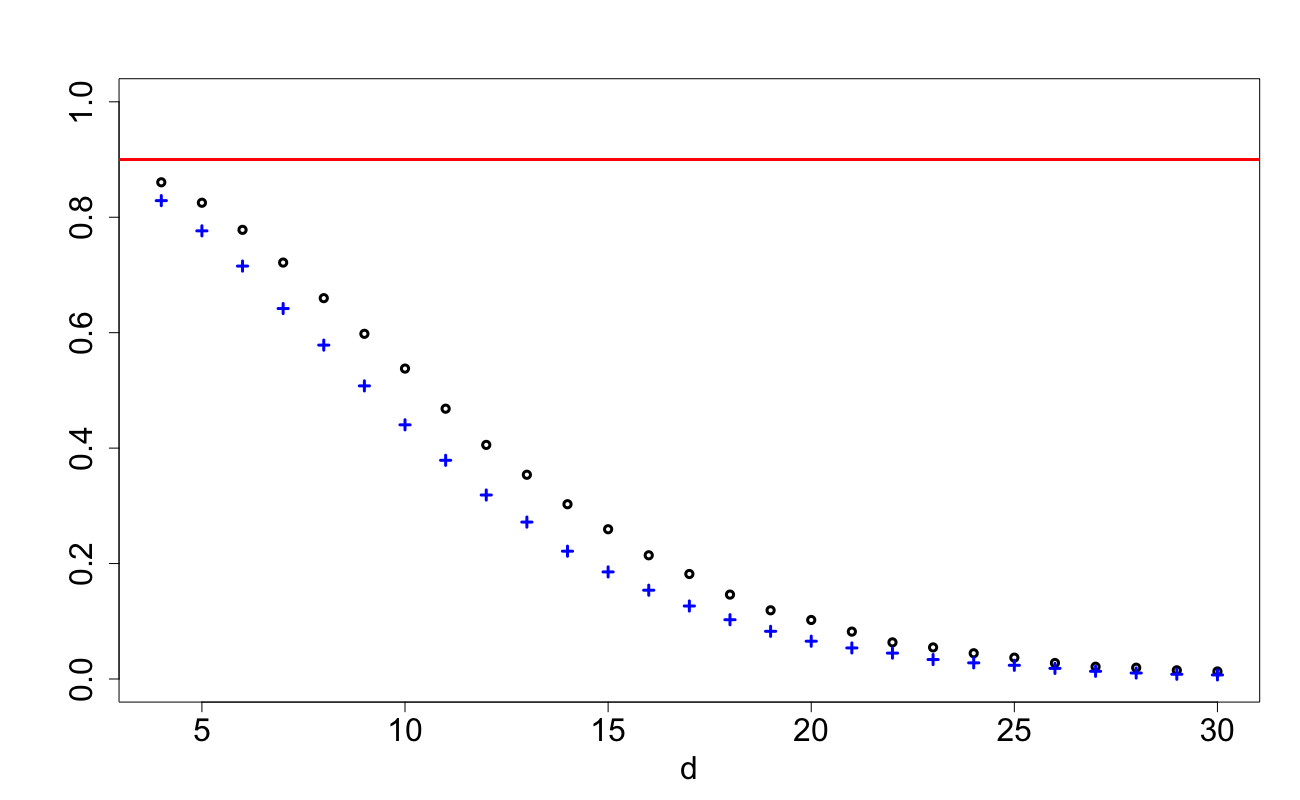}
  \caption{Covering using asymptotic \\radius: $n=1000,10000$. }
  \label{key_figure1}
\end{minipage}%
\begin{minipage}{.5\textwidth}
  \centering
  \includegraphics[width=1\linewidth]{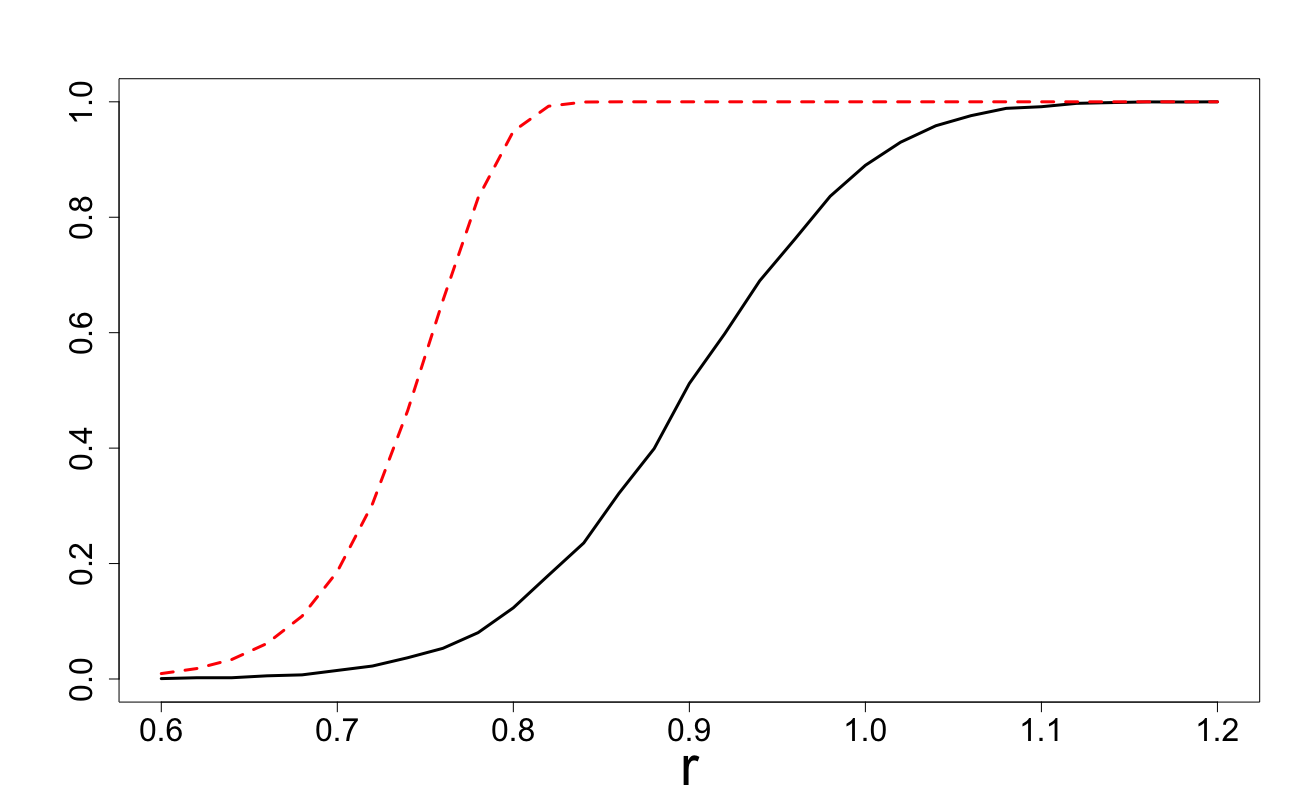}
  \caption{Weak covering and its asymptotic approximation: $d=20, n=10000$ }
    \label{key_figure2}
\end{minipage}
\end{figure}

For the popular scenario of $\X=[0,1]^d$, in Figure~\ref{key_figure1} for $n=1,000$ (blue plusses) and $n=10,000$ (black circles) we depict $F (r_{asy},X_n)$ as a function of $d$, where $r_{asy}$ is the asymptotic radius obtained from \eqref{eq:Zador1} to achieve 0.9 covering. We see that very quickly and for $n$ that would be deemed large, the true weak covering is significantly less than 0.9 and quickly tends to zero in $d$. The big difference between the asymptotic and finite regime is further illustrated in Figure~\ref{key_figure2}. Using a solid black line we depict $F (r,X_n)$ as a function of $r$ with $d=20$ and $n=10000$. In this figure, the dashed red line is the approximation obtained from the asymptotic result \eqref{eq:Zador1}, that is, the approximation $F(r,X_n) \approx F_d(n^{1/d}V_d^{1/d}r)$.

\begin{figure}[!h]
\centering
\begin{minipage}{.5\textwidth}
  \centering
  \includegraphics[width=1\linewidth]{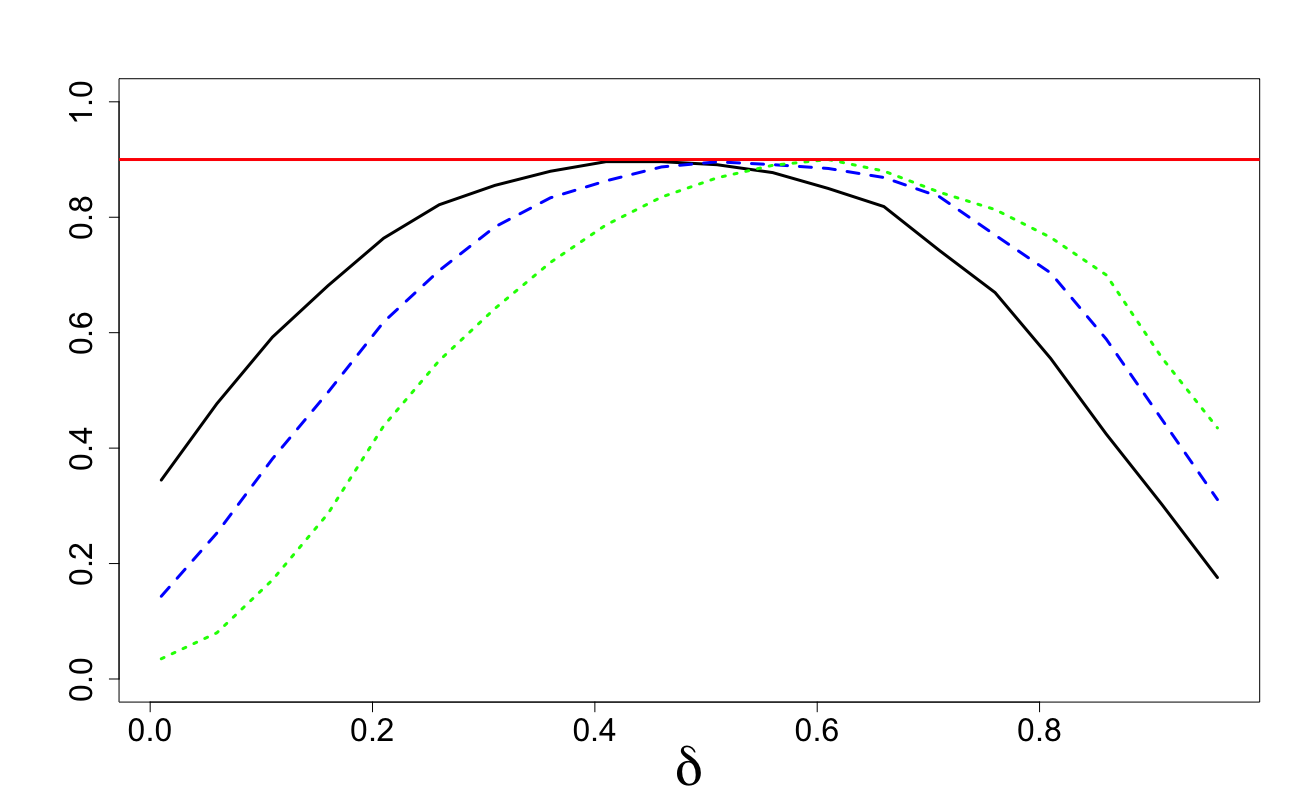}
  \caption{Weak covering of $\X=[-1,1]^d$:\\ $d=50, n=1000,10000,100000$. } \label{Main_pic2}
  \label{Main_pic1}
\end{minipage}%
\begin{minipage}{.5\textwidth}
  \centering
  \includegraphics[width=1\linewidth]{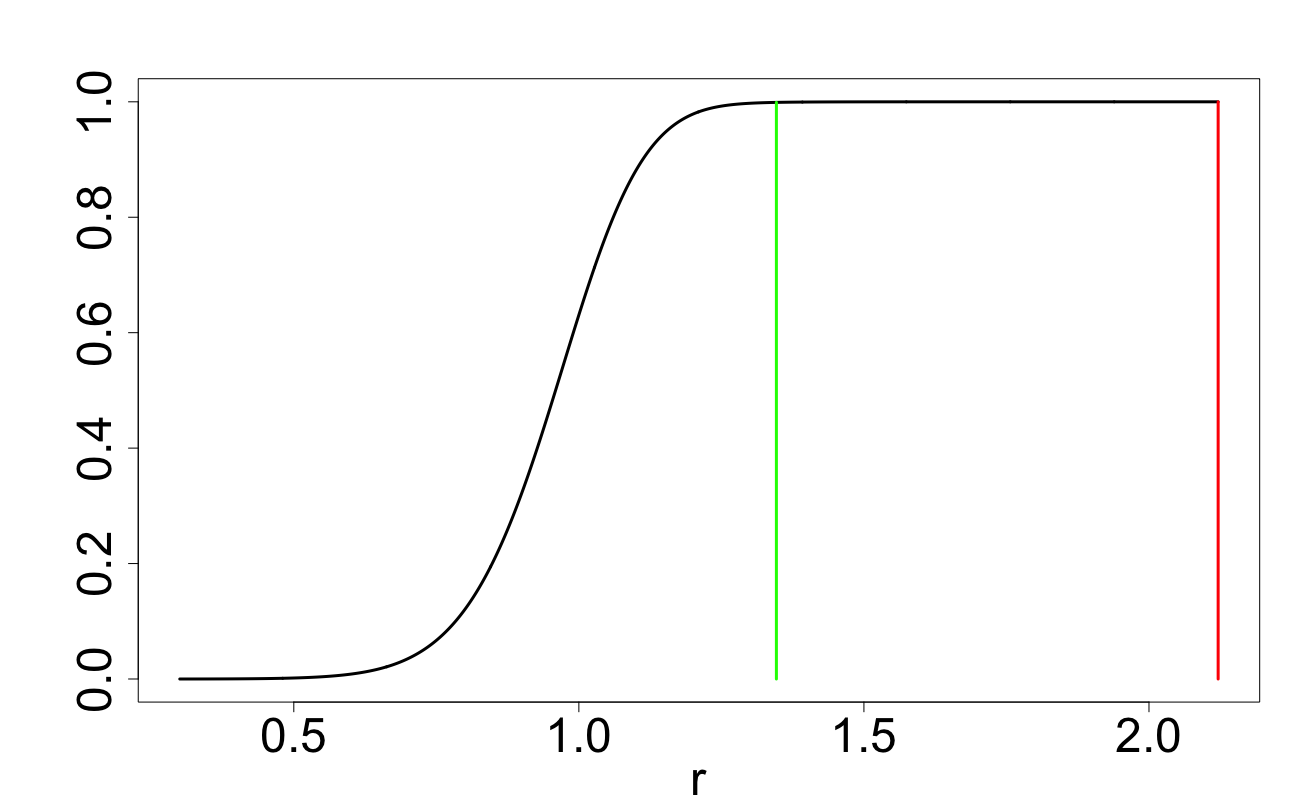}
  \caption{$F(r,X_{n})$ for $2^{d-1}$ design with $r_{0.999}$ and $r_{1}$: $d=10$}
    \label{Main_pic2}
\end{minipage}
\end{figure}

For $\X=[-1,1]^d$, the weak covering properties of a number of random and deterministic designs was studied in \cite{us} and \cite{second_paper}. The second phenomena of weak covering of high-dimensional sets is the so-called `$\delta$-effect'.
The $\delta$-effect is a phenomenon in high dimensions that results in the recommendation to restrict your design $X_n$ to within the cube $[-\delta,\delta]^d$ with $0<\delta<1$, instead of $[-1,1]^d$. An example of the $\delta$-effect is demonstrated in Figure~\ref{Main_pic1}. Here, for $d=50$, $X_n$ is a sample of $n$ i.i.d. uniform random vectors from the $\delta$ cube $[-\delta,\delta]^d$. The $x$-axis corresponds to the value of $\delta$ and the $y$-axis corresponds to the proportion of $\X$ covered. For $n=1000$ (solid black), $10000$ (dashed blue), $100000$ (dotted green), the value of $r$ has been selected to ensure the covering at optimal $\delta$ is 0.9. From this picture it is clear that in high dimensions, even with a large number of design points, it is beneficial for weak covering to choose your design in a $\delta$-cube with $\delta <<1$. The case of $\delta=1$ leads to poor weak covering.

 In Figure~\ref{Main_pic2} we demonstrate another phenomenon of covering the high-dimensional cube $\X=[-1,1]^d$ (along with the $\delta$-effect). The third phenomena highlights the difficulty and excessive nature of requiring the full covering of $\X$. Motivated by the numerical results of \cite{us}, in \cite{noonan2022efficient} a theoretical investigation into a $2^{d-1} $ design of maximum resolution concentrated at the
  points
  $(\pm 1/2, \ldots, \pm 1/2) \in \mathbb{R}^d$ was performed. For this design, the c.d.f. $F (\cdot,X_n)$ is shown with a black line and we also indicate the location of the $r_{1}={\rm CR}(X_n)$ and $r_{0.999}$ by vertical red and green line respectively. In this figure, we  take $d=10$ and therefore $n=512$.
It is very easy to analytically compute the covering radius (for any $d>2$):
  ${\rm CR}(X_n)=\sqrt{d+8}/2 $ ; for $d=10$ this gives
  ${\rm CR}(X_n) \simeq 2.1213 $.
The value of  $r_{0.999}$ satisfies $r_{0.999}(X_n) \simeq 1.3465 $. This value has been computed using very accurate approximations developed in \cite{noonan2022efficient}; we claim  3 correct decimal places in  $r_{0.999}$. This figure demonstrates Theorem 3.4 in \cite{noonan2022efficient} which states for any $0<\gamma<1$ and for this special contruction $X_n$, ${r_{1-\gamma}/r_1 \rightarrow 1/\sqrt{3}}$ as $d\rightarrow \infty$. So in large dimensions, covering, for example, $99.99\%$ of $\X$ can be achieved with a radius approximately 0.577 times smaller than the full covering radius.

\subsection*{Conclusion}

The problem of covering of a given compact set $\X \subset \mathbb{R}^d$ is formulated and several issues  occurring in the case when $d$ is large are discussed. It is demonstrated that the common covering schemes like uniform random or pseudo-random sampling in $\X$  make poor covering but there are some  ways of improving such coverings. 
The results discussed can make implications for devising exploration strategies in wide variety of global optimization algorithms in  high dimensions.

\paragraph*{Related entries from within the Encyclopedia of Optimization:}

 Random search for global optimization; Random search methods; Convergence of global random search algorithms.

\bibliographystyle{unsrt}


\bibliography{large_dimension_enc_1}

\begin{thebibliography}{10}

\bibitem{schaback2006kernel}
R.~Schaback and H.~Wendland.
\newblock Kernel techniques: from machine learning to meshless methods.
\newblock {\em Acta numerica}, 15:543--639, 2006.

\bibitem{wendland2004scattered}
H.~Wendland.
\newblock {\em Scattered data approximation}, volume~17.
\newblock Cambridge university press, 2004.

\bibitem{Niederreiter}
H.~Niederreiter.
\newblock {\em Random number generation and quasi-{M}onte {C}arlo methods}.
\newblock SIAM, Philadelphia, PA, 1992.

\bibitem{pronzato2012design}
L.~Pronzato and W.~M{\"u}ller.
\newblock Design of computer experiments: space filling and beyond.
\newblock {\em Statistics and Computing}, 22(3):681--701, 2012.

\bibitem{santner2003design}
T.~Santner, B.~Williams, W.~Notz, and B.~Williams.
\newblock {\em The design and analysis of computer experiments}, volume~1.
\newblock Springer, 2003.

\bibitem{penrose2021random}
M.~Penrose.
\newblock Random euclidean coverage from within.
\newblock {\em arXiv preprint arXiv:2101.06306}, 2021.

\bibitem{Sukh1}
A.~Sukharev.
\newblock Optimal strategies of search for an extremum.
\newblock {\em USSR Computational Mathematics and Math. Physics},
  11(4):910--924, 1971.

\bibitem{sukharev2012minimax}
A.~Sukharev.
\newblock {\em Minimax models in the theory of numerical methods}.
\newblock Springer, 1992.

\bibitem{vzilinskas2013worst}
A.~{\v{Z}}ilinskas.
\newblock On the worst-case optimal multi-objective global optimization.
\newblock {\em Optimization Letters}, 7(8):1921--1928, 2013.

\bibitem{pages1998space}
G.~Pag{\`e}s.
\newblock A space quantization method for numerical integration.
\newblock {\em Journal of computational and applied mathematics}, 89(1):1--38,
  1998.

\bibitem{du1999centroidal}
Q.~Du, V.~Faber, and M.~Gunzburger.
\newblock Centroidal voronoi tessellations: Applications and algorithms.
\newblock {\em SIAM review}, 41(4):637--676, 1999.

\bibitem{borodachov2019discrete}
S.~Borodachov, D.~Hardin, and E.~Saff.
\newblock {\em Discrete energy on rectifiable sets}.
\newblock Springer, 2019.

\bibitem{brass2005research}
P.~Brass, W.~Moser, and J.~Pach.
\newblock {\em Research problems in discrete geometry}, volume~18.
\newblock Springer, 2005.

\bibitem{boroczky2004finite}
K.~B{\"o}r{\"o}czky~Jr.
\newblock {\em Finite packing and covering}, volume 154.
\newblock Cambridge University Press, 2004.

\bibitem{kuperbergball}
G~Kuperberg and W~Kuperberg.
\newblock Ball packings and coverings with respect to the unit cube (in
  preparation).

\bibitem{pronzato2019measures}
L.~Pronzato and A.~Zhigljavsky.
\newblock Measures minimizing regularized dispersion.
\newblock {\em J. of Scientific Computing}, 78(3):1550--1570, 2019.

\bibitem{zhigljavsky2021bayesian}
A.~Zhigljavsky and A.~{\v{Z}}ilinskas.
\newblock {\em Bayesian and high-dimensional global optimization}.
\newblock Springer, 2021.

\bibitem{pronzato2020bayesian}
L.~Pronzato and A.~Zhigljavsky.
\newblock Bayesian quadrature, energy minimization, and space-filling design.
\newblock {\em SIAM/ASA Journal on Uncertainty Quantification}, 8(3):959--1011,
  2020.

\bibitem{graf2007foundations}
S.~Graf and H.~Luschgy.
\newblock {\em Foundations of quantization for probability distributions}.
\newblock Springer, 2007.

\bibitem{us}
J.~Noonan and A.~Zhigljavsky.
\newblock Covering of high-dimensional cubes and quantization.
\newblock {\em SN Operations Research Forum}, 1(3):1--32, 2020.

\bibitem{second_paper}
J.~Noonan and A.~Zhigljavsky.
\newblock Non-lattice covering and quantization of high dimensional sets.
\newblock In {\em Black Box Optimization, Machine Learning, and No-Free Lunch
  Theorems}, pages 273--318. Springer, 2021.

\bibitem{noonan2022efficient}
J.~Noonan and A.~Zhigljavsky.
\newblock Efficient quantisation and weak covering of high dimensional cubes.
\newblock {\em Discrete \& Computational Geometry}, pages 1--26, 2022.

\end{thebibliography}

\newpage

\end{document}